\documentclass[11pt,twoside,a4paper]{article}

\usepackage{amsmath}
\usepackage{amsthm}
\usepackage{changepage}
\usepackage{geometry}
\usepackage{txfonts}
\usepackage{cite} 
\usepackage{listings}
\usepackage{hyperref}
\usepackage{amssymb}
\usepackage{mathrsfs}
\author{}%}Liu Yong, Liu Ziyu}
\hypersetup{
	colorlinks=true,
	linkcolor=blue,
	anchorcolor=blue,
	citecolor=blue}
\title{ \textbf{Asymptotic stability for non-equicontinuous
Markov semigroups}}
\newtheorem{theorem}{Theorem}[section]

\numberwithin{equation}{section}
\theoremstyle{definition}
\newtheorem{remark}[theorem]{Remark}
\newtheorem{example}[theorem]{Example}
\newtheorem{definition}[theorem]{Definition}

\begin{document}

 {\author{ \bf{Fuzhou Gong$^a$,\ Yong Liu$^b$,\ Yuan Liu$^a$,\   Ziyu Liu$^b$} \\ 
 \small{\em  {\bf$^a$} Academy of Mathematics and Systems Science, Chinese Academy of Sciences, Beijing 100190, China}\\
\small{\em  {\bf$^b$}   LMAM, School of Mathematical Science, Peking University, 100871, Beijing, China}\\
\small{\em Email:\;fzgong@amt.ac.cn(F. Gong);\ liuyong@math.pku.edu.cn (Yong Liu);}\\
\small{\em liuyuan@amss.ac.cn(Yuan Liu);\ liuziyu@math.pku.edu.cn(Z. Liu)}
}}

	\maketitle

	\makeatother

   \begin{abstract}

We prove that the asymptotic stability, also known as the weak mixing, is equivalent to a lower bound condition together with the eventual continuity. The latter is a form of weak regularity for Markov-Feller semigroups that generalizes the e-property. Additionally, we provide an example of an asymptotically stable Markov semigroup with essential randomness that does not satisfy the e-property.
\end{abstract}
 
\maketitle

\section{Introduction}\label{Sec 1}
    This paper mainly centers on the  \emph{asymptotic stability} (also known as the \textit{weak mixing}, see \cite[Section 11.2]{DPZ2014}) of Markov-Feller semigroups. Specifically, we adopt the notion of \emph{asymptotic equicontinuity condition}, introduced by Jaroszewska \cite{J2013}, or the \emph{eventual continuity}, by the first and third author of the present paper \cite{GL2015} (see Definition \ref{Def eventual continuity}), to  work as a necessary condition for the asymptotic stability.  In fact, these two notions are formulated almost simultaneously, and mathematically equivalent. We adopt the name of eventual continuity for simplicity, which depicts the feature that a uniquely ergodic semigroup may behave sensitively in initial data, and is evidently weaker than the  \emph{e-property} (see Definition \ref{Def e-property}), i.e. equicontinuity. The eventual continuity is a reasonable and essential tool to deal with non-equicontinuous semigroups. \par

    In 2006, Lasota and Szarek developed the lower bound technique for \emph{e-processes} (i.e., Markov processes with the e-property) to formulate a criterion for the existence of an invariant measure \cite{LS2006,S2006} on a non-locally compact metric space. Since then, the e-property has turned out to play a valuable role in the study of ergodicity.  For example, Szarek and Worm use the lower bound condition to formulate a criterion for the asymptotic stability in \cite[Corollary 5.4]{SW2012} from the Markov-Feller semigroups with the e-property. In the opposite direction, in 2017, it is shown that any asymptotically stable Markov-Feller semigroup satisfies the e-property if the support of its unique invariant measure has an interior, see \cite{HSZ2017}. More recently, further investigations into this direction have been conducted by Kukulski and Wojew\'{o}dka-\'{S}ci\c{a}\.{z}ko, see \cite{KW2021,KW2024} for a detailed account. Let us mention that since the e-property requires some uniform estimates in a neighbourhood of each point, it does not hold in general cases, or even if it does, it can be cumbersome to verify, see the discussion in \cite{NJG2017,KS2018}. On the other hand, a recent paper \cite{LL2024} by the second and last author characterizes the relation between  the e-property and the eventual continuity. Moreover, it is proved that these two properties are equivalent restricted on the support of each ergodic measure.  \par 
    In the present paper, the lower bound condition is applied to provide an equivalent  characterization between the asymptotic stability and the eventual continuity (see Theorem \ref{Thm 1}). We also construct a non-trivial example with randomness that satisfies asymptotic stability, and from which eventual continuity follows, but it does not satisfy the e-property. More precisely, this example is a variant inspired by a class of the iterated function systems considered in  \cite{BKS2014}. See more about the iterated function systems in \cite{HMS2005,HMS2006}. It is worth mentioning that, to the best of our knowledge, the known examples of the non-euqicontinuous Markov-Feller semigroups are  deterministic dynamical systems, see \cite{GL2015,HSZ2017}. \par 
 
    This paper is organized as follows. Section \ref{Sec 2} gives some basic facts and notations of Markov-Feller semigroups. The main result is presented in Section \ref{Sec 3}.  In Section \ref{Sec 4}, we provide an iterated function system example that is asymptotically stable, but does not satisfy the e-property.

    \section{Preliminaries}\label{Sec 2}

    In this section, we introduce some basic facts and  several definitions of Markov-Feller semigroups. Let $(\mathcal{X},\rho)$ be a Polish space. By $\mathcal{B}(\mathcal{X})$ we denote the $\sigma$-field of all its Borel subsets. By $B(x, r)$ we denote the open ball in $(\mathcal{X},\rho)$ of radius $r$, centred at $x\in\mathcal{X}$.  By $B_b(\mathcal{X})$, we denote the vector space of all bounded real-valued Borel measurable functions on $\mathcal{X}$, and by $C_b(\mathcal{X})$ all bounded real-valued continuous functions, both equipped with the supremum norm $|\cdot|_{\infty}$. By $L_b(\mathcal{X})$ we denote the subspace of $C_b(\mathcal{X})$ of all bounded Lipschitz functions.  By $\mathcal{M}_1(\mathcal{X})$, we denote the family of all probability measures on $\mathcal{X}$. For $\mu\in\mathcal{M}_1(\mathcal{X})$, its \emph{support} is the set
	$\text{supp }\mu=\{x\in \mathcal{X}:\mu(B(x,\epsilon))>0 \text{ for every } \epsilon>0\}$. Moreover, for $f\in B_b(\mathcal{X})$ and $\mu\in\mathcal{M}_1(\mathcal{X})$, we use the notation $\langle f,\mu\rangle=\int_{\mathcal{X}}f(x)\mu(dx).$	\par
	Let $T$ be the index set, $\mathbb{R}_+$ or $\mathbb{N}_+$. $\{P_t\}_{t\in T}$ is a Markov-Feller semigroup on $\mathcal{X},$ if $\{P_t\}_{t\in T}$ is a Markov semigroup and $P_t(C_b(\mathcal{X}))\subset C_b(\mathcal{X})$ for $t\in T$. Let $\{P^*_t\}_{t\in T}$ be the dual semigroup defined on $\mathcal{M}_1(\mathcal{X})$ by the formula $P^*_t\mu(B):=\int_{\mathcal{X}}P_t\mathbf{1}_Bd\mu$ for $B\in\mathcal{B}(\mathcal{X}).$ For ease of notation, we simply rewrite $P^*_t$ to $P_t$. Recall that $\mu\in\mathcal{M}_1(\mathcal{X})$ is invariant for the semigroup $\{P_t\}_{t\in T}$ if $P_t\mu=\mu$ for all $t\in T$. We now introduce the notion of asymptotic stability, which is also called weak mixing, see e.g. \cite[Section 11.2]{DPZ2014}.
    \begin{definition}\label{Def asymptotic stability}
	   $\{P_t\}_{t\in T}$ is asymptotically stable if there exists a unique invariant measure $\mu\in\mathcal{M}_1(\mathcal{X})$, and  $P_t\nu$ converges weakly to $\mu$  as $t\rightarrow\infty$ for every $\nu\in\mathcal{M}_1(\mathcal{X})$.
    \end{definition}

   Finally, let us recall the following two types of regularities of Markov semigroups, namely, the e-property and the eventual continuity. 
    
   \begin{definition}\label{Def e-property}
	    $\{P_t\}_{t\in T}$ satisfies the e-property, if for any $z\in \mathcal{X}$ and $f\in L_b(\mathcal{X})$,
	\begin{equation*}
    \lim\limits_{x\rightarrow z}\sup\limits_{t\geq 0}|P_tf(x)-P_tf(z)|=0.
	\end{equation*}
	\end{definition}

    \begin{definition}\label{Def eventual continuity}
	    $\{P_t\}_{t\in T}$ is eventually continuous, if for any $z\in\mathcal{X}$ and $f\in L_b(\mathcal{X})$,
	\begin{equation*}
	\lim\limits_{x\rightarrow z}\limsup\limits_{t\rightarrow\infty}|P_tf(x)-P_tf(z)|=0.
	\end{equation*}
	\end{definition}

\section{Main results}\label{Sec 3}

 Throughout this paper, we assume that $\{P_t\}_{t\in T}$ is a Markov-Feller semigroup. One interesting observation gives that the eventual continuity is necessary for the asymptotic stability, which can be proved by merely comparing the definitions between these two notions.   Moreover, for eventually continuous semigroups, the asymptotic stability can be derived from some lower bound condition on probability transition.  Our principal result is a new criterion of asymptotic stability for the eventual continuous semigroups as follows.

	\begin{theorem}\label{Thm 1} 
		The following two statements are equivalent:
			\begin{itemize}
				\item[$(a)$] $\{P_t\}_{t\in T}$ is asymptotically stable with a unique invariant measure $\mu$.
								
				\item[$(b)$] $\{P_t\}_{t\in T}$ is eventually continuous, and there exists some $z\in\mathcal{X}$ such that for any $\epsilon>0$, \begin{equation}\label{eq E1}
				\inf\limits_{x\in\mathcal{X}}\liminf\limits_{t\rightarrow\infty}P_t(x,B(z,\epsilon))>0.
				\end{equation}
        
		\end{itemize}		
	\end{theorem}
    Before turning to the proof of Theorem \ref{Thm 1}, let us first make some comments on the results.

   \begin{remark}
       The implication $(b)\Rightarrow (a)$ extends the criterion for asymptotic stability in \cite[Corollary 5.4]{SW2012} from the e-processes to eventually continuous semigroups. It is worth mentioning that, in our case, we actually derive an equivalent relation between the asymptotic stability and the eventual continuity.
   \end{remark}
     \begin{proof}[Proof of Theorem~{\upshape\ref{Thm 1}}]
 We first prove that $(a)\Rightarrow(b)$. In view of Definition \ref{Def asymptotic stability}, for any $x,z\in\mathcal{X}$ and $f\in L_b(\mathcal{X})$, one gets
     \begin{equation*}
	|P_tf(x)-P_tf(z)|\leq |P_tf(x)-\langle f,\mu\rangle| + |P_tf(z)-\langle f,\mu\rangle|\rightarrow 0\quad\text{ as }t\rightarrow\infty,
	\end{equation*}
     which clearly implies the eventual continuity.\par 
     In order to prove the lower bound condition, let us fix any $z\in\text{supp }\mu$ and $\epsilon>0$. Then for any $x\in\mathcal{X}$, it follows that
     \begin{equation*}
         \lim\limits_{t\rightarrow\infty}P_t(x,B(z,\epsilon))=\mu(B(z,\epsilon))>0.
     \end{equation*}
     Thus the desired inequality (\ref{eq E1}) follows from the fact that the positive lower bound $\mu(B(z,\epsilon))$ is independent of $x\in\mathcal{X}$.\par

    Next we show that $(b)\Rightarrow(a)$, which consists of two parts: {\bf Part 1} on the asymptotic behavior of semigroups and {\bf Part 2} on the existence of an invariant measure.\par 
    
    {\bf Part 1}.  First we show that $\lim\limits_{t\rightarrow\infty}|P_tf(x)-P_tf(y)|=0$ for all $x,y\in \mathcal{X}$ and all $f\in L_b(\mathcal{X}).$ Assume, contrary to our claim, there exists  some  $x_1,x_2\in \mathcal{X}$, $f\in L_b(\mathcal{X})$ and $\epsilon>0$ such that $\limsup\limits_{t\rightarrow\infty}|P_tf(x_1)-P_tf(x_2)|\geq 3\epsilon$. As $\{P_t\}_{t\in T}$ is eventually continuous, there is $\delta>0$ such that
		\begin{equation*}
		\limsup\limits_{t\rightarrow\infty}|P_tf(x)-P_tf(z)|<\epsilon/2\quad\text{for any }x\in B(z,\delta). 
		\end{equation*}
		For such $\delta>0$, condition (\ref{eq E1}) implies that there exists a positive number $\alpha=\alpha(\delta) \in (0,\frac12)$ such that $\liminf\limits_{t\rightarrow\infty}P_t(x,B(z,\delta))>\alpha$ for all $x\in\mathcal{X}$.	Then Fatou's lemma gives, for any $\nu\in\mathcal{M}_1(\mathcal{X})$,
		\begin{equation}\label{eq alpha}
		\liminf\limits_{t\rightarrow\infty}P_t\nu(B(z,\delta))\geq\int_{\mathcal{X}}\liminf\limits_{t\rightarrow\infty}P_t(y,B(z,\delta))\nu(dy)>\alpha.
		\end{equation}\par 
		Let $k\geq 1$ be such that $2(1-\alpha)^k|f|_{\infty}<\epsilon.$ By induction we are going to define four sequences of measures $\{\nu_i^{x_1}\}_{i=1}^k,\{\mu_i^{x_1}\}_{i=1}^k,\{\nu_i^{x_2}\}_{i=1}^k,\{\mu_i^{x_2}\}_{i=1}^k,$ and a sequence of positive numbers $\{t_{i}\}_{i=1}^k$ in the following way.  In view of inequality (\ref{eq alpha}), let us choose $\nu=\delta_{x_j}\in\mathcal{M}_1(\mathcal{X})$, $j=1,2$, then there exists some 
		$t_1=t_1(x_1,x_2)>0$  sufficiently large such that 	$P_{t}(x_j,B(z,\delta))>\alpha$	for all $t\geq t_1$, $j=1,2$.
		Set \begin{center}
			$\nu_1^{x_j}(\cdot) = \dfrac{P_{t_1}\delta_{x_j}(\cdot \cap B(z,\delta))}{P_{t_1}\delta_{x_j}(B(z,\delta))},\quad$
			$\mu_1^{x_j}(\cdot) = \dfrac{1}{1-\alpha}(P_{t_1}\delta_{x_j}(\cdot)-\alpha\nu_1^{x_j}(\cdot)),\quad j=1,2.$
		\end{center}
		Assume that we have done it for $i = 1,\dots , l,$ for some $l < k.$ Now using (\ref{eq alpha}) again, let $t_{l+1}\geq t_l$ be sufficiently large such that
		$P_{t}\mu_l^{x_j}(B(z,\delta))>\alpha$ for all $t\geq t_{l+1}$, $j=1,2$. 
		Set 
		\begin{center}
			$\nu_{l+1}^{x_j}(\cdot) = \dfrac{P_{t_{l+1}}\mu_l^{x_j}(\cdot \cap B(z,\delta))}{	P_{t_{l+1}}\mu_l^{x_j}(B(z,\delta))},\quad$
			$\mu_{l+1}^{x_j}(\cdot) = \dfrac{1}{1-\alpha}(P_{t_{l+1}}\mu_l^{x_j}(\cdot)-\alpha\nu_{l+1}^{x_j}(\cdot)),\quad j=1,2.$
		\end{center}
		Then it follows that
		\begin{equation*}
		\begin{aligned}
		P_{t_1+\dots+t_k}\delta_{x_j}(\cdot)&=\alpha P_{t_2+\dots+t_k}\nu_1^{x_j}(\cdot)+\alpha(1-\alpha) P_{t_3+\dots+t_k}\nu_2^{x_j}(\cdot)+\dots+\\
         &  \quad+\alpha(1-\alpha)^{k-1}\nu_k^{x_j}(\cdot)+(1-\alpha)^k \mu_k^{x_j}(\cdot),
		\end{aligned}
		\end{equation*}
		where $\text{supp }{{\nu}_i^{x_j}}\subset B(z,\delta),\,i=1,\dots,k,\,j=1,2.$ Thus  from the Fatou's lemma and eventual continuity,  we have
		\begin{equation*}
		\begin{aligned}
	 \limsup\limits_{t\rightarrow\infty}|\langle P_{t}f, \nu_i^{x_1} \rangle-\langle P_{t}f, \nu_i^{x_2} \rangle| &\leq \limsup\limits_{t\rightarrow\infty}|\langle P_{t}f-P_{t}f(z), \nu_i^{x_1} \rangle| +\\
    & \quad\,\limsup\limits_{t\rightarrow\infty}|\langle P_{t}f-P_{t}f(z), \nu_i^{x_2} \rangle|\\
    &\leq\epsilon/2+\epsilon/2=\epsilon.
		\end{aligned}
		\end{equation*}\par 
		Finally,  using the measure decomposition gives
		\begin{equation*}
		\begin{aligned}
		3\epsilon&\leq\limsup\limits_{t\rightarrow\infty}|P_tf(x_1)-P_tf(x_2)|\\
		&=\limsup\limits_{t \to \infty}|\langle P_tf,P_{t_1+\cdots+t_k} \delta_{x_1}\rangle-\langle P_tf,P_{t_1+\cdots+t_k} \delta_{x_2}\rangle| \\
		&\leq \alpha\limsup\limits_{t \to \infty}|\langle P_{t}f, \nu_1^{x_1} \rangle-\langle P_{t}f, \nu_1^{x_2} \rangle|+\cdots+\\
            &\quad+\alpha(1-\alpha)^k\limsup\limits_{t \to \infty}|\langle P_{t}f, \nu_k^{x_1}\rangle-\langle P_{t}f, \nu_k^{x_2} \rangle|+\epsilon\\
		&\leq (\alpha+\cdots+\alpha(1-\alpha)^k)\epsilon+\epsilon< 2\epsilon,
		\end{aligned}
		\end{equation*}
		which is impossible.
    
    {\bf Part 2}. It remains to establish the existence of an invariant measure.   Invoking a standard Kryloy-Bogolyubov procedure, the existence of an invariant measure will be implied by the tightness of the sequence $\{P_t(z,\cdot)\}_{t\in T}$. We argue by contradiction.    Assuming, on the contrary, that $\{P_t(z,\cdot)\}_{t\in T}$ is not tight, we first show that there exists a positive number $\epsilon$, a sequence of compact sets $\{K_i:i\geq 1\}$, and an increasing sequence of real numbers $\{t_i:i\geq 1\}$ tending to infinity, such that 
    \begin{equation}\label{Ki-def1}
        P_{t_i}(z,K_i)\geq \epsilon,\quad\text{ for any }i\geq 1,
    \end{equation}
    \begin{equation}\label{Ki-def2}
        K_i^{\epsilon/2}\cap K_j^{\epsilon/2}=\emptyset\quad\text{ for any }i\neq j,
    \end{equation}
    where we use the notation $A^{\epsilon}:=\{x\in\mathcal{X}:\inf_{y\in A} \rho(x,y)<\epsilon\}$ for $A\subset \mathcal{X}$.\par 

    As we have assumed that $\{P_t(z,\cdot)\}_{t\in T}$ is not tight, thus there exists some positive number $\epsilon$ such that for any compact set $K\subset\mathcal{X}$ and any $t_0\geq 0$, there exists $t=t(K)\geq t_0$ such that
    \begin{equation}\label{Ki-def3}
        P_t(z,K^{\epsilon})<1-\epsilon.
    \end{equation}
    Let us fix any compact set $F_1\subset\mathcal{X}$ together with $t_1\geq 0$ such that  $P_{t_1}(z,F_1^{\epsilon})<1-\epsilon$.  Then  we can choose a compact set $K_1\subset\mathcal{X}\setminus F_1^\epsilon$ such that
    \begin{equation*}
        P_{t_1}(z,K_1^{\epsilon})>\epsilon.
    \end{equation*}
   We  inductively construct the sequence $\{K_i:i\geq 1\}$ and $\{t_i:i\geq 1\}$ together with an auxiliary sequence of compact sets $\{F_i:i\geq 1\}$ satisfying (\ref{Ki-def1}), (\ref{Ki-def2}). Assume that we have done it for $i = 1,\dots,n$ for some $n\geq 1$. Let $F_{n+1}=F_n\cup K_n$. By (\ref{Ki-def3}), we may find $t_{n+1}\geq t_{n}+1$ such that
   \begin{equation*}
        P_{t_{n+1}}(z,F_{n+1}^{\epsilon})<1-\epsilon.
    \end{equation*}
    Then we again let $K_{n+1}\subset \mathcal{X}\setminus F_{n+1}^\epsilon$ be a compact set such that 
    \begin{equation*}
        P_{t_{n+1}}(z,K_{n+1}^{\epsilon})>\epsilon,
    \end{equation*}
    which completes the construction. In particular, this construction ensures (\ref{Ki-def1}), (\ref{Ki-def2}) and that $\{t_i:i\geq 1\}$ increases to infinity. \par   
    We now define $f_i\in L_b(\mathcal{X})$ by
   \begin{equation*}
    f_i(y):=\rho(y,(K_i^{\epsilon/4})^c)/(\rho(y,(K_i^{\epsilon/4})^c)+\rho(y,K_i)),
    \end{equation*}
    satisfying that
    \begin{equation}\label{eq 1}	
			\mathbf{1}_{K_{i}}\leq f_i\leq\mathbf{1}_{K_{i}^{\epsilon/4}}\quad\text{and}\quad\text{Lip}(f_i)\leq 4/\epsilon \ \ \ \ \textrm{for any}\ i\in\mathbb{N}.	
    \end{equation}
    
    We claim that there exists a subsequence of integers $\{j_k: k\geq 0\}$ and some point $y_0\in \mathcal{X}$ such that for  $f:=\sum_{l\geq 0}f_{j_l}\in L_b(\mathcal{X})$
    \begin{equation}\label{eq 1`}
        	\limsup\limits_{k \to \infty}P_{t_{j_k}}f(y_0)\leq\epsilon/2.
    \end{equation} \par 
    Admitting this claim (\ref{eq 1`}) for the moment, it follows from the definition of $f_{j_k}$ that
         \begin{equation*}
    	P_{t_{j_k}}f(z) \geq P_{t_{j_k}} f_{j_k}(z) \geq P_{t_{j_k}}(z,K_{j_k}) \geq \epsilon > \epsilon/2\ \ \ \ \textrm{for all}\ k\geq 0,
          \end{equation*}	
    which contradicts the asymptotic result in {\bf Part 1}. This completes the proof of tightness.\par 
	
	Following some arguments in \cite[Theorem 1.6]{GL2015}, we prove the claim (\ref{eq 1`}) as follows.\par 
	    $\mathbf{Step\;1}.$  Denote $B_0=B(z,r).$ By (\ref{eq E1}), let
		$\alpha_0:=\liminf\limits_{i\rightarrow\infty} P_{t_i}(z,B_0)>0$.
		Lemma 2.2 in \cite{GL2015} yields that there exists $j_0$ sufficiently large such that
    \begin{equation*}
        \begin{aligned}            \epsilon\alpha_0/16\geq\liminf\limits_{i\rightarrow\infty}P_{t_i}(z,K_{j_0}^{\epsilon/4})&=\liminf\limits_{i\rightarrow\infty}\int_{\mathcal{X}} P_{t_{j_0}}(y,K_{j_0}^{\epsilon/4})P_{t_i}(z,dy)\\
        &\geq\liminf\limits_{i\rightarrow\infty}\int_{B_0} P_{t_{j_0}}f_{j_0}(y)P_{t_i}(z,dy).
        \end{aligned}
    \end{equation*}

		Due to the Feller property, define an open subset $A_0=\{y\in B_0:P_{t_{j_0}}f_{j_0}(y)<\epsilon/8\}$.
		It follows that
		\begin{equation*}
	 \epsilon\alpha_0/16\geq\liminf\limits_{i\rightarrow\infty}\int_{B_0-A_0} P_{t_{j_0}}f_{j_0}(y)P_{t_i}(z,dy) \geq \liminf\limits_{i\rightarrow\infty}P_{t_i}(z,B_0-A_0) \cdot \epsilon/8,
		\end{equation*}
		which implies $\liminf\limits_{i\rightarrow\infty}P_{t_i}(z,B_0-A_0)\leq\alpha_0/2$, and thus $\limsup\limits_{i\rightarrow\infty}P_{t_i}(z,A_0)\geq\alpha_0/2$.
		Hence $A_0$ is nonempty, which contains a ball $B_1$ of radius less than $r/2$  such that $P_{t_{j_0}}f_{j_0}(y)\leq\epsilon/8$ for all $y\in \overline{B_1}$  and $P_s(z,B_1)>0$ for some $s>0$. By the Feller property, there exists $\delta>0$ such that for any  $z'\in B(z,\delta),$ $P_s(z',B_1)\geq\frac{1}{2}P_s(z,B_1)>0$.
		From  Condition (\ref{eq E1}), we derive that
         \begin{equation*}
        \begin{aligned} 
     \alpha_1:=\liminf\limits_{t\rightarrow\infty}P_t(z,B_1)&\geq \liminf\limits_{t\rightarrow\infty} \int_{B(z,\delta)} P_s(z',B_1) P_t(z,dz')\\
     &\geq\frac{1}{2}P_s(z,B_1)\cdot\liminf\limits_{t\rightarrow\infty}P_t(z,B(z,\delta))>0.
         \end{aligned}
     \end{equation*}\par 
		Inductively, for each $k\geq 0,$ we can determine  $B_{k+1}\subset B_k$ of radius less than $r/2^{k+1}$ and  increasing $j_k\in\mathbb{N}$ such that $P_{t_{j_k}}f_{j_k}(y)\leq\epsilon/8$ for all $ y\in \overline{B_{k+1}}$. Therefore, we can find a common $y_0\in\bigcap_{k\geq 0}\overline{B_{k}}$ such that
		\begin{equation}\label{eq 2`}
		P_{t_{j_k}}f_{j_k}(y_0)\leq\epsilon/8,\quad\forall k\geq 0.
		\end{equation}\par 
		
	 $\mathbf{Step\;2}.$
	     Noting that all $K_i^{\epsilon/2}$ are disjoint mutually, there exists some $u_0$ sufficiently large such that $P_{t_{j_0}}\sum_{l\geq u_0}f_{j_l}(y_0)\leq\epsilon/8$.
		Combining with (\ref{eq 2`}) yields $P_{t_{j_0}}(f_{t_{j_0}}+\sum_{l\geq u_0}f({j_l}))(y_0)\leq\epsilon/4$. Similarly, there exists $u_1>u_0$ such that $P_{t_{j_{u_0}}}(f_{t_{j_{u_0}}}+\sum_{l\geq u_1}f({j_l}))(y_0)\leq\epsilon/4$.
		By induction, we obtain a subsequence $\{j_0, j_{u_0}, j_{u_1}, \ldots\} \subset\{j_l\}$, still denoted by $\{j_l\}$ for simplicity, such that
	    \begin{equation}\label{eq 4`}
		P_{t_{j_k}}\sum_{l\geq k}f_{j_l}(y_0)\leq\epsilon/4,\quad\forall k\geq 0.
		\end{equation}\par
		
		On the other hand, by \cite[Lemma 2.3]{GL2015} and (\ref{eq 2`}), applying diagonal arguments, passing to a subsequence if necessary, we further obtain 
		\begin{equation}\label{eq 3`}
		\limsup\limits_{k\rightarrow\infty}P_{t_{j_{k}}}f_{j_{l}}(y_0)\leq\epsilon/(8\cdot 2^l),\quad\forall l\geq 0.
		\end{equation}\par
		By (\ref{eq 3`}), there exists some $v_0$ sufficiently large such that 
		$P_{t_{j_k}}f_{j_0}(y_0)\leq\epsilon/8$ for all $k\geq v_0$. Similarly, there exists $v_1>v_0$ such that $P_{t_{j_{k}}}(f_{j_{0}}+f_{j_{v_0}})(y_0)\leq(1+2^{-1})\epsilon/8$ for all $k\geq v_1$.	 By induction, we derive a subsequence $\{j_0,j_{v_0},j_{v_1}, \ldots\}\subset\{j_l\}$, still denoted by $\{j_l\}$, satisfying that
		\begin{equation}\label{eq 5`}
		P_{t_{j_k}}\sum_{0\leq l<k}f_{j_l}(y_0)\leq \sum_{0\leq  l <k} 2^{-l} \epsilon/8 \leq\epsilon/4.
		\end{equation}
		Combining (\ref{eq 5`}) with  (\ref{eq 4`}), it yields 
		\begin{equation*}
		P_{t_{j_k}}\sum_{l\geq 0}f_{j_l}(y_0)\leq\epsilon/2,\quad\forall k\geq 0,
		\end{equation*}	
    which completes the proof of our claim (\ref{eq 1`}).
	\end{proof}

	\begin{example}    We adopt the idea in \cite[Example 2.1]{HSZ2017} to construct a non-equicontinuous  Markov chain in continuous-time. Let $\mathcal{X}=\{0\}\cup\{n^{-1}:n\geq 2\}\cup\{n:n\geq 2\}.$  Define a continuous-time Markov chain $\{\Phi_t\}_{t\geq 0}$ with transition rate matrix $\{q_{ij}\}_{i,j\in \mathcal{X}},$ where $q_{n,n}=q_{n^{-1},n^{-1}}=-n,\;q_{n,0}=q_{n^{-1},n}=n$, for $n\geq 2,$ and $q_{i,j}=0$ otherwise. The transition probability $\{p_{ij}(t)\}_{i,j\in \mathcal{X}}$ is defined as follows: $p_{0,0}(t)=1$; for $n\geq 2$, 
	\begin{align*}
	&p_{n^{-1},n^{-1}}(t)=p_{n,n}(t)=e^{-t/n},  &p_{n^{-1},n}(t)=\frac{t}{n}e^{-t/n}, \\ &p_{n^{-1},0}(t)=1-p_{n^{-1},n^{-1}}(t)-p_{n^{-1},n}(t), &p_{n,0}(t)=1-e^{-t/n}.
	\end{align*}
	It can be  checked that $\{P_t\}_{t\geq 0}$ generated by $\{\Phi_t\}_{t\geq 0}$ is a Markov-Feller semigroup. Further, since $\Phi_t$ converges to $0$ almost surely as $t\rightarrow\infty$, $\{P_t\}_{t\geq 0}$ is asymptotically stable, hence $\{P_t\}_{t\geq 0}$ is eventually continuous by Theorem~{\upshape\ref{Thm 1}}. However, the e-property fails at $0$, since we can take $f(x)=x\land 1 \in C_b(\mathcal{X})$ and then
		\begin{equation*}
    	P_{n}f(n^{-1})-P_{n}f(0)\geq e^{-1},\;\forall\, n\geq 2.
		\end{equation*}
    \end{example} 

\section{A non-equicontinuous Markov-Feller semigroup}\label{Sec 4}

 In this part, we construct an asymptotically stable Markov-Feller semigroup with essential randomness, which does not satisfy the e-property (see Example \ref{Ex 1}). Specifically, we modify the  iterated function systems (IFS) introduced in \cite{BKS2014}.  The randomness comes from that, at each step,  the dynamical system jumps to different positions, obeying some place-dependent probabilistic distribution.    \par 
	We first introduce some notations. Let $(\mathcal{X},\rho)$ be a Polish space, $I=\{1,\dots,N\}$ be a finite set, $w_i:\mathcal{X}\rightarrow\mathcal{X},\;i\in I$ be transformations on $\mathcal{X}$, and $\{S(t)\}_{t\geq 0}$ be a continuous flow on $\mathcal{X}$. Let $(p_1(x),\dots,p_N(x))$ be a probabilistic vector  denoted by $p(x)$, that is, for each $x\in \mathcal{X}$, $p_i(x)\geq 0$ for $i\in I$ and $\sum_{i\in I}p_i(x)=1$. Let $\{\tau_n\}_{n\geq 1}$  be a sequence of random variables with $\tau_0=0,$ and $\triangle\tau_n=\tau_n-\tau_{n-1},n\geq 1$, are i.i.d. with density $\lambda e^{-\lambda t},\;\lambda>0$.\par 
	Now we define the $\mathcal{X}$-valued Markov chain $\Phi=\{\Phi^x(t):x\in \mathcal{X}\}_{t\geq 0}$ in the following way. 
	
	\begin{itemize}
	    \item[$\mathbf{Step\;1.}$] Let $x\in\mathcal{X}$. Denote $\Phi_0^x:=x$.
	    \item[$\mathbf{Step\;2.}$] Let $\xi_1=S(\tau_1)(x).$  We   have a random index $i_1\in I$ with probability $p(\xi_1)$, i.e. $
		\mathbb{P}(i_1=k)=p_{k}(\xi_1)$, for $k\in I$. 
		Set $\Phi_1^x:=w_{i_1}(\xi_1).$
	    \item[$\mathbf{Step\;3.}$]  Recursively, assume that $\Phi_1^x,\dots,\Phi_{n-1}^x,n\geq 2$ are given. We define $\xi_n=S(\triangle\tau_n)(\Phi_{n-1}^x).$ Further, we randomly choose $i_n\in I$ with  probability $p(\xi_n)$, and let $\Phi_n^x=w_{i_n}(\xi_n).$
	    \item[$\mathbf{Step\;4.}$]  Finally, set  $
		\Phi^x(t):=S(t-\tau_n)(\Phi_n^x)\quad\text{for  }\tau_n\leq t<\tau_{n+1},\;n\geq 0.$
	\end{itemize}
	
	Now we can define semigroup $\{P_t\}_{t\geq 0}$ by $P_tf(x):=\mathbb{E}f(\Phi^x(t))\text{ for } f\in B_b(\mathcal{X})$.

	Given the above setting, we provide the following example with essential randomness which is not an e-process, but is eventually continuous.

	\begin{example}\label{Ex 1}
		Let $\mathcal{X}=\mathbb{R}_+.\;$ Let $S(t)= \text{Id}_{\mathcal{X}}$ for $t\geq 0,$ and $ w_1(x)=0,$ $w_2(x)=x,$ $w_3(x)=x^{-1}\mathbf{1}_{\{x\neq 0\}}$. Let $p_i\in C(\mathcal{X},[0,1]),$ $i=1,2,3$, be such that
		\begin{equation*}
		(p_1(x),p_2(x),p_3(x))=\begin{cases}
		(\frac{x}{2},1-x,\frac{x}{2}),&0\leq x< \frac{2}{3},\\
		(\frac{1}{3},\frac{1}{3},\frac{1}{3}),&\frac{2}{3}\leq x\leq \frac{3}{2},\\
		(\frac{1}{2x},1-x^{-1},\frac{1}{2x}),&x>\frac{3}{2}.
		\end{cases}
		\end{equation*}
		By construction, any trajectory starting from $x\neq 0$ runs over $\{0,x,x^{-1}\}$.	Noting that $0$ is the unique absorbing state for this process,  $\Phi^x_t$ converges to $0$ almost surely for all initial states $x\in \mathcal{X}$ as $t\rightarrow\infty.$ 	Therefore $\{P_t\}_{t\geq 0}$  is asymptotically stable and  eventually continuous by Theorem~{\upshape\ref{Thm 1}}. However, it does not satisfy the e-property.\par 
		In fact, choosing $x_n=1/n$, we have that
		\begin{equation*}
		\begin{aligned}
		\mathbb{P}(\Phi_t^{x_n}=n)&\geq\sum_{k=1}^{\infty}\sum_{l=0}^{k-1}p_2(x_n)^lp_3(x_n)p_2(x_n^{-1})^{k-l-1}\mathbb{P}(\tau_k\leq t<\tau_{k+1})\\
		&=\frac{x_n}{2}\sum_{k=1}^{\infty}(1-x_n)^{k-1}\frac{(\lambda t)^k}{(k-1)!}e^{-\lambda t}\\
		&=\frac{1}{2}x_n\lambda te^{-x_n\lambda t}.
		\end{aligned}
		\end{equation*}
		Now, letting $f(x)=x\wedge 1\in C_b(\mathcal{X})$ and $t=n$, it follows that, 
		\begin{equation*}
		P_nf(x_n)\geq \mathbb{P}(\Phi_n^{x_n}=n)f(n)\geq \frac{1}{2}\lambda e^{-\lambda }>0=P_nf(0),
		\end{equation*}
		which conflicts the e-property.

	\end{example}

		{\small \section*{Acknowledgement}  Fuzhou Gong and Yuan Liu are supported by National Key R \& D Program of China (No. 2020YFA0712700) and Key Laboratory of Random Complex Structures and Data Sciences, Mathematics and Systems Science, Chinese Academy of Sciences (No. 2008DP173182); Yong Liu is supported by CNNSF (No. 11731009, No. 12231002) and Center for Statistical Science, PKU; Yuan Liu is supported by  (CNNSF 11688101).}

\bibliographystyle{amsplain}

\begin{thebibliography}{10}



\bibitem{BKS2014}
H. Bessaih, R. Kapica, and T. Szarek, \textit{ Criterion on stability for {M}arkov processes applied to a model with jumps}, Semigroup Forum \textbf{88} (2014), no. 1, 76--92. 

\bibitem{NJG2017}
N. Glatt-Holtz, J. C. Mattingly and G. Richards,  \textit{On unique ergodicity in nonlinear stochastic partial differential equations}, J. Stat. Phys. \textbf{166} (2017), no. 3-4, 618 -- 649. 

\bibitem{DPZ2014}
G. Da Prato, J. Zabczyk,  \textit{Stochastic equations in infinite demensions}, Second edition, Cambridge University Press, Cambridge (2014). 


\bibitem{GL2015}
F. Gong and Y. Liu, \textit{Ergodicity and asymptotic stability of {F}eller semigroups on {P}olish metric spaces}, Sci. China Math. \textbf{58} (2015), no. 6, 1235--1250. 


\bibitem{HMS2005}
K. Horbacz, J. Myjak and T. Szarek, \textit{On stability of some general random dynamical system},  J. Stat. Phys. \textbf{119}  (2005), no. 1-2, 35-60. 

\bibitem{HMS2006}
K. Horbacz, J. Myjak and T. Szarek, \textit{Stability of random dynamical systems on Banach spaces},  Positivity \textbf{10}  (2006), no. 3, 517–538. 

\bibitem{HSZ2017}
S. C. Hille, T. Szarek and M. A. Ziemla\'nska, \textit{Equicontinuous families of {M}arkov operators in view of asymptotic stability},  C. R. Math. Acad. Sci. Paris \textbf{355}  (2017), no. 12, 1247--1251.


\bibitem{KW2021}
R. Kukulski and H. Wojew\'{o}dka-\'{S}ci\c{a}\.{z}ko, \textit{The e-property of asymptotically stable Markov-Feller operators},  Colloq. Math.  \textbf{165} (2021), no. 2, 269--283. 

\bibitem{KW2024}
R. Kukulski and H. Wojew\'{o}dka-\'{S}ci\c{a}\.{z}ko, \textit{TThe e-property of asymptotically stable Markov semigroups},  Results Math.  \textbf{79} (2024), no. 3, Paper No. 112, 22 pp. 



\bibitem{J2013}
J. Jaroszewska, \textit{On asymptotic equicontinuity of {M}arkov transition functions},  Statist. Probab. Lett.  \textbf{83} (2013), no. 3, 943--951. 



\bibitem{KS2018}
A. Kulik and M. Scheutzow, \textit{Generalized couplings and convergence of transition probabilities}, 
  Probab. Theory Related Fields  \textbf{171} (2018), no. 1-2, 333 -- 376. 

\bibitem{LL2024}
Y. Liu and Z. Liu, \textit{Relation between the eventual continuity and the e-property for Markov-Feller semigroups},  Acta Math. Appl. Sin. Engl. Ser.  \textbf{40} (2024), no. 1, 1 -- 16.

\bibitem{LS2006}
A. Lasota and T. Szarek, \textit{Lower bound technique in the theory of a stochastic differential equation},  J. Differential Equations \textbf{231} (2006), no. 2, 513 -- 533. 

\bibitem{S2006}
T. Szarek, \textit{{F}eller processes on nonlocally compact spaces},  Ann. Probab.  \textbf{34} (2006), no. 5, 1849 -- 1863. 

\bibitem{SW2012}
T. Szarek and D. T. H. Worm, \textit{Ergodic measures of {M}arkov semigroups with the e-property},  Ergodic Theory Dynam. Systems \textbf{143} (2012), no. 3, 1117--1135. 


\end{thebibliography}

\end{document}